\documentclass[12pt, a4paper, reqno]{amsart}
\usepackage{mathrsfs}
\usepackage[active]{srcltx}
\usepackage{mathrsfs}
\usepackage{amssymb}

\usepackage[left=2.5cm, right=2.5cm, top=2.50cm, bottom=2.5cm]{geometry}

\newtheorem{proposition}{Proposition}[section]
\newtheorem{theorem}{Theorem}[section]
\newtheorem{lemma}{Lemma}[section]

\newtheorem{remark}{Remark}[section]

\numberwithin{equation}{section}

\newcommand{\Urt}{{\bf U}}

\usepackage{color}
\usepackage[unicode]{hyperref}
\hypersetup{
	colorlinks = true,%
	citecolor = [rgb]{0.0,0.0,0.9},
	filecolor=black,%
	linkcolor = [rgb]{0.65,0.0,0.0},%
	anchorcolor = red,
	pagecolor = red,
	urlcolor= [rgb]{0.65,0.0,0.0}
}

\begin{document}
	
	\title[]{%Sharp Heisenberg-Pauli-Weyl 
Nodal solutions for the fractional Yamabe problem on Heisenberg groups}
	
	%    Information for first author
	\author{Alexandru Krist\'aly}
	%    Address of record for the research reported here
\noindent	\address{\noindent Department of Economics, Babe\c s-Bolyai University,   400591 Cluj-Napoca,
		Romania \& Institute of Applied Mathematics, \'Obuda
		University,
		  1034 Budapest, Hungary}
	\email{alexandrukristaly@yahoo.com; kristaly.alexandru@nik.uni-obuda.hu}
	%    \thanks will become a 1st page footnote.
%	\thanks{Research supported by OTKA, project no. 124827.  
		%CNCS-UEFISCDI, project no. PN-II-ID-PCE-2011-3-0241. %, "Symmetries
		%in elliptic problems: Euclidean and non-Euclidean techniques"
		%and by the J\'anos Bolyai Research Scholarship of the Hungarian
		%Academy of Sciences. 
		%The present work was initiated during the
		%author's visit at the Institut des Hautes \'Etudes Scientifiques
		%(IH\'ES), Bures-sur-Yvette, France.
%	}

	%    General info

	
	\begin{abstract}
	We prove that  the fractional Yamabe equation $\mathcal L_\gamma u=|u|^\frac{4\gamma}{Q-2\gamma}u$ on the Heisenberg group $\mathbb H^n$ has  $[\frac{n+1}{2}]$ sequences of nodal (sign-changing) weak solutions whose elements have mutually different nodal properties, where  $\mathcal L_\gamma$ denotes the 
		CR fractional sub-Laplacian operator on $\mathbb H^n$,    $Q=2n+2$ is the homogeneous dimension of $\mathbb H^n$, and $\gamma\in \bigcup_{k=1}^n[k,\frac{kQ}{Q-1})$.
		Our argument is  variational, based on a Ding-type conformal pulling-back transformation of the original  problem into a problem on the CR sphere $S^{2n+1}$ combined with a suitable Hebey-Vaugon-type compactness result and   group-theoretical constructions for special subgroups of the unitary group ${\bf U}(n+1).$
	\end{abstract}
	
		\dedicatory{Dedicated to 
		Professor Patrizia Pucci on the occasion of her 65th birthday}
	
	\maketitle

	\section{Introduction}
	
After the seminal paper of Caffarelli and Silvestre \cite{Caf-Silv}, considerable efforts have been made concerning the study of elliptic problems involving the
fractional Laplace operator both in Euclidean and non-Euclidean settings. As expected, the Euclidean framework is much more developed; although many results concerning the pure Laplace operator can be transposed to the fractional setting in $\mathbb R^n$, there are also subtle differences which require a deep understanding of certain nonlinear phenomena, see e.g. Cabr\'e and Sire \cite{Cabre}, Caffarelli \cite{Caf}, Caffarelli,  Salsa and Silvestre \cite{Caf-Salsa-Silv}, Di Nezza,  Palatucci and Valdinoci \cite{Valdinoci}, and references therein.

By exploring analytical and spectral theoretical arguments, important contributions have been obtained recently within the CR setting concerning the fractional Laplace operator with various applications in sub-elliptic PDEs, see Branson, Fontana and Morpurgo \cite{BFM}, Frank and Lieb \cite{FL}, and Frank, del Mar Gonz\'alez, Monticelli and Tan \cite{Frank-Adv-Math}. In particular, in the latter papers sharp Sobolev and Moser-Trudinger inequalities are established  on the Heisenberg group $\mathbb H^n$, $n\geq 1,$ the simplest non-trivial CR structure. 

In the present paper we shall consider the \textit{fractional Yamabe problem} on the Heisenberg group $\mathbb H^n$, namely, 
\begin{equation*}
\left\{
\begin{array}{lll}
\mathcal L_\gamma u=|u|^\frac{4\gamma}{Q-2\gamma}u & \mbox{on} & \mathbb H^n, \\
u\in  D^{\gamma}({\mathbb{H}^n}).
\end{array}%
\right.\eqno{{(\textbf{FYH})_\gamma}}
\end{equation*}
Hereafter, $Q:=Q_n=2n+2$ is the  homogeneous dimension of $\mathbb H^n$,  $\gamma>0$ is a parameter specified later,  %$0<\gamma<\frac{Q}{2}=n+1$,  
$\mathcal L_\gamma$ denotes the CR fractional sub-Laplacian operator on $\mathbb H^n,$ and the functional space $D^{\gamma}({\mathbb{H}^n})$ contains real-valued functions from $L^{\frac{2Q}{Q-2\gamma}}(\mathbb H^n)$ whose energy associated to the CR fractional sub-Laplacian operator $\mathcal L_\gamma$ is finite; see \S \ref{subsection2-2} for details.    

Due to the  recent paper of Frank, del Mar Gonz\'alez, Monticelli and Tan \cite{Frank-Adv-Math}, we know the existence of positive solutions of  ${(\textbf{FYH})_\gamma}$  for  $\gamma\in (0,Q/2)$, having the form  
\begin{equation}\label{positive}
u(z,t)=c_0\left((1+|z|^2)^2+t^2\right)^\frac{2\gamma -Q}{4},\ \ (z,t)\in \mathbb H^n,
\end{equation}
for some $c_0>0,$ allowing any left translation and dilation.
%, where $c_0=\left({\Gamma\left(\frac{n+1-\gamma}{2}\right)}/{\Gamma\left(\frac{n+1+\gamma}{2}\right)}\right)^\frac{2(Q-2\gamma)}{Q\gamma-Q-2\gamma^2-2\gamma}.$ 
In the special case $\gamma=1,$ when $\mathcal L_1=\mathcal L$ is the usual sub-Laplacian operator on $\mathbb H^n,$ the existence and uniqueness (up to left translation and dilation) of positive solutions of the form  (\ref{positive}) for problem ${(\textbf{FYH})_1}$ have been established by Jerison and Lee \cite{Jerison-Lee-1, Jerison-Lee-2}; see also Garofalo and Vassilev \cite{Garofalo} for generic Heisenberg-type groups (e.g. Iwasawa groups). 

%\medskip

Our main result guarantees sign-changing solutions for the fractional Yamabe problem ${(\textbf{FYH})_\gamma}$ as follows: 

\begin{theorem}\label{fotetel}
	Let
	%$0<\gamma<\frac{Q}{2}=n+1$ with $[\gamma(1-\frac{1}{Q})]+1=[\gamma]$. 
	$\gamma\in \bigcup_{k=1}^n[k,\frac{kQ}{Q-1})$, where $Q=2n+2.$ 
	Then problem {\rm ${(\textbf{FYH})_\gamma}$} admits  at least $[\frac{n+1}{2}]$ sequences of sign-changing weak solutions whose elements have mutually different nodal properties. $($Hereafter, $[r]$ denotes the integer part of $r\geq 0$.$)$
\end{theorem}

%Some remarks are in order. % concerning Theorem \ref{fotetel}. 

Before commenting on Theorem \ref{fotetel} we recall that similar results are well known in the Euclidean setting; indeed, Bartsch, Schneider and Weth \cite{BSW} proved the existence of infinitely many sign-changing weak solutions for the polyharmonic problem 
\begin{equation*}
	\left\{
	\begin{array}{lll}
		(-\Delta)^mu=|u|^\frac{4m}{N-2m}u & \mbox{in} & \mathbb R^N, \\
		u\in  \mathcal D^{m,2}({\mathbb{R}^N}),
	\end{array}%
	\right.\eqno{{(\textbf{P})_m}}
\end{equation*}
where $N>2m$, $m\in \mathbb N$, and  $\mathcal D^{m,2}({\mathbb{R}^N})$ denotes the usual higher order Sobolev space over $\mathbb R^N.$ In fact, their proof is based on Ding's original idea, see \cite{Ding}, who considered the case $m=1$, by pulling back the variational problem ${(\textbf{P})_m}$ to the standard sphere $S^{N}$ by stereographic projection. In this manner, by exploring certain properties of suitable subgroups of  %$\textsf{O}(N_1)\times \textsf{O}(N_2)$ of 
the orthogonal group ${\bf O}(N+1)$,  the authors are able to obtain compactness by exploring 
a suitable Sobolev embedding result of Hebey and Vaugon \cite{HV} which is indispensable in the application of the symmetric mountain pass theorem. 

We  notice that sign-changing solutions are already guaranteed to the usual CR-Yamabe problem ${(\textbf{FYH})_1}$  by Maalaoui and  Martino \cite{Martino-0}, and Maalaoui,   Martino and Tralli \cite{Martino-1} by exploring  Ding's approach; their results are direct consequences of Theorem \ref{fotetel} for  $\gamma=1$.

Coming back to  Theorem \ref{fotetel}, we shall mimic Ding's original idea as well, emphasizing that our \textit{CR fractional} setting requires a more delicate analysis than either the  polyharmonic setting in the Euclidean case (see \cite{BSW}) or the usual CR framework, i.e., when $\gamma=1$ (see \cite{Martino-0, Martino-1}). In the sequel we sketch our strategy. As expected, we first  consider the fractional Yamabe problem on the CR sphere $S^{2n+1}$, i.e., 
\begin{equation*}
	\left\{
	\begin{array}{lll}
		\mathcal A_\gamma U=|U|^\frac{4\gamma}{Q-2\gamma}U & \mbox{on} & S^{2n+1}, \\
		U\in H^\gamma(S^{2n+1}),
	\end{array}%
	\right.\eqno{{(\textbf{FYS})_\gamma}}
\end{equation*}
where the intertwining operator $\mathcal A_\gamma$ and Sobolev space $H^\gamma(S^{2n+1})$ are introduced in \S \ref{subsection2-2}. By using the Cayley transform between the Heisenberg group $\mathbb H^n$ and the CR sphere $S^{2n+1}$ we prove that there is an explicit correspondence between the weak solutions of {\rm ${(\textbf{FYH})_\gamma}$} and {\rm ${(\textbf{FYS})_\gamma}$}, respectively, see Proposition \ref{weak-equivalent} (and Remark \ref{weak-2} for an alternative proof). Being in the critical case, the energy functional associated to problem ${(\textbf{FYS})_\gamma}$ does not satisfy the usual Palais-Smale condition due to the lack of compactness of the embedding $H^\gamma(S^{2n+1})\hookrightarrow L^\frac{2Q}{Q-2\gamma}(S^{2n+1})$. In order to regain some compactness, we establish a CR fractional version of the Ding-Hebey-Vaugon compactness result on the CR sphere $S^{2n+1}$, see Proposition \ref{proposition-compact}. In fact,  subgroups of the unitary group ${\bf U}(n+1)$ having the form $G={\bf U}(n_1)\times...\times {\bf U}(n_k)$  with $n_1+...+n_k=n+1$
%$G={\bf U}(l)\times {\bf U}(m)$ (with $l,m\geq $ and $l+m=n+1$) 
will imply the compactness of the embedding of $G$-invariant functions of $H^\gamma(S^{2n+1})$ into $L^\frac{2Q}{Q-2\gamma}(S^{2n+1})$. Here, we shall explore 
the compactness result of Maalaoui and  Martino \cite{Martino-0} combined with an iterative argument of  Aubin \cite{Aubin} and the technical assumption   $\gamma\in \bigcup_{k=1}^n[k,\frac{kQ}{Q-1})$; some comments on the necessity of the latter  assumption are formulated in Remark \ref{remark-kell}.
%; although the latter assumption is crucial in our compactness argument,  we believe it can be removed by using another approach. 
%In fact, the (sharp) CR fractional Sobolev inequality on $S^{2n+1}$ by Frank and Lieb \cite{FL} and specific subgroups of the unitary group ${\bf U}(n+1)$ having the form $G={\bf U}(k)\times {\bf U}(m)$ (with $k+m=n+1$) will imply the compactness of the embedding of $G$-invariant functions of $H^\gamma(S^{2n+1})$ into $L^\frac{2Q}{Q-2\gamma}(S^{2n+1})$. Here, the most fastidious part of the compactness proof is the  fine splitting properties of spherical/zonal harmonics on $S^{2n+1}$  with respect to the underlying group $G$, where some representations of $G$-invariant functions are deeply explored from  Vilenkin and Klimyk \cite{Vilenkin}. 
Now, having such a compactness, the fountain theorem and the principle of symmetric criticality applied to the energy functional associated to ${(\textbf{FYS})_\gamma}$ will guarantee the existence of a whole sequence of $G$-invariant  weak solutions for  {\rm ${(\textbf{FYS})_\gamma}$}, so for {\rm ${(\textbf{FYH})_\gamma}$}. The number of  $[\frac{n+1}{2}]$ sequences of sign-changing weak solutions for  {\rm ${(\textbf{FYH})_\gamma}$} with mutually different nodal properties will follow by  careful choices of the subgroups $G={\bf U}(n_1)\times...\times {\bf U}(n_k)$ of the unitary group ${\bf U}(n+1)$ with $n_1+...+n_k=n+1$, see Proposition \ref{proposition-Rubik}. 

\medskip

\textit{Plan of the paper.} In Section \ref{section-2} we recall those notions and results that are indispensable to present our argument (e.g. basic facts about Heisenberg groups, the Cayley transform, spherical/zonal harmonics on $S^{2n+1}$, fractional Sobolev spaces on $S^{2n+1}$ and $\mathbb H^n$). Section \ref{section-3} is devoted to the proof of Theorem \ref{fotetel}; in \S \ref{section-3-1} we prove the equivalence between the weak solutions of problems {\rm ${(\textbf{FYS})_\gamma}$} and {\rm ${(\textbf{FYH})_\gamma}$}; in \S \ref{section-3-2} we establish the compactness result on the CR fractional setting for $S^{2n+1}$; in \S \ref{section-3-3}  we treat the group-theoretical aspects of our problem concerning the choice of the subgroups $G={\bf U}(n_1)\times...\times {\bf U}(n_k)$ of the unitary group ${\bf U}(n+1)$ which is needed to produce $[\frac{n+1}{2}]$ sequences of sign-changing weak solutions for {\rm ${(\textbf{FYH})_\gamma}$} with different nodal properties. Finally, in \S \ref{section-3-4} we assemble the aforementioned pieces in order to conclude the proof of Theorem \ref{fotetel}.

%; to do so, we shall roughly explore the idea of solving the Rubik cube, see Balogh and Krist\'aly \cite{BK}. 

%\medskip

%\textit{Plan of the paper.} ...

	\section{Preliminaries}\label{section-2}
	
In order the paper to be self-contained,  we recall in this section some basic notions from \cite{BFM}, \cite{Folland-PAMS},  \cite{FL},  \cite{Frank-Adv-Math} and \cite{Vilenkin} which are indispensable in the sequel.
	
	\subsection{Heisenberg groups}

	An element in the Heisenberg group $\mathbb H^n$ is denoted by  $(x, y, t)  $, where $x
	= (x_1, \ldots, x_n) \in \mathbb R^n$, $y = (y_1, \ldots, y_n)
	\in \mathbb R^n$, $t \in \mathbb R,$ and we identify the pair $(x, y)$ with
	$z \in \mathbb C^n$ having coordinates $z_j = x_j + i y_j$ for all $j = 1, \ldots, n$.
	The correspondence with its Lie algebra via the exponential
	coordinates induces the  group law
	$$(z, t) \star (z', t') = \left( z + z', t + t' + 2{\rm Im}\  z \cdot \overline{z'} \right), \quad \forall (z, t),\ (z', t') \in \mathbb C^n \times \mathbb R,$$
	where ${\rm Im}$ denotes the imaginary part of a complex number and
	$ z \cdot \overline{z'} =\sum\limits_{j = 1}^n z_j
	\overline{z_j'}$ is the Hermitian inner product. In these
	coordinates the neutral element of $\mathbb H^n$ is $0_{\mathbb
		H^n}=(0_{\mathbb C^n},0)$ and the inverse  $(z, t)^{-1}$ of  the element $(z, t)$ is
	$(-z, -t)$. Note that $(x, y, t) = (z, t)$ forms a
	real coordinate system for $\mathbb H^n$ and the system of vector
	fields given as differential operators
	$$X_j = \frac{\partial}{\partial{x_j}} + 2 y_j \frac{\partial}{\partial t}, \quad Y_j = \frac{\partial}{\partial{y_j}} - 2 x_j \frac{\partial}{\partial t}, \quad  j \in \{ 1, \ldots n\},\ \ T = \frac{\partial}{\partial t},$$
	forms a basis of the left-invariant vector fields on $\mathbb H^n.$
	The vectors $X_j,Y_j, j \in \{1,...,n\}$ form the basis of the
	horizontal bundle.  Let $$N(z,t)=(|z|^4+t^2)^\frac{1}{4}$$ be the homogeneous gauge norm on $\mathbb H^n$ and $d_{KC}:\mathbb H^n\times \mathbb H^n\to \mathbb R$ be the Kor\'anyi-Cygan metric given by
	$$d_{KC}((z,t),(z',t'))=N((z', t')^{-1} \star (z, t))=(|z-z'|^4+(t-t'-2 {\rm Im}\ z \cdot \overline{z'} )^2)^\frac{1}{4}.$$
	The Lebesgue measure of $\mathbb R^{2n+1}$ will be
	the Haar measure on $\mathbb H^n$ (uniquely defined up to a positive
	multiplicative constant). 

\subsection{Cayley transform}
	Let $$S^{2n+1}=\{\zeta=(\zeta_1,...,\zeta_{n+1})\in \mathbb C^{n+1}:\langle \zeta, \overline \zeta\rangle=1\}$$ be the unit sphere in $\mathbb C^{n+1}$, where $\langle\zeta, \overline \eta \rangle=\sum_{j=1}^{n+1}\zeta_j\overline {\eta_j}.$ The distance function on $S^{2n+1}$ is given by $$d_S(\zeta,\eta)=\sqrt{2|1-\langle\zeta, \overline \eta \rangle|},\ \ \ \zeta,\eta\in S^{2n+1}.$$

	The Cayley transform 
	$\mathcal C:\mathbb H^n\to  S^{2n+1}\setminus \{(0,...,0,-1)\}$ is defined by 
	$$\mathcal C(z,t)=\left(\frac{2z}{1+|z|^2+it},\frac{1-|z|^2-it}{1+|z|^2+it}\right),$$
whose Jacobian determinant 
 is given by 
 $${\rm Jac}_{\mathcal C}(z,t)=\frac{2^{2n+1}}{((1+|z|^2)^2+t^2)^{n+1}},\ \ \ (z,t)\in \mathbb H^n.$$	
 Accordingly, for any integrable function $f:S^{2n+1}\to \mathbb R$, we have 
\begin{equation}\label{valtozo-csere}
\int_{S^{2n+1}}f(\eta){\rm d}\eta=\int_{\mathbb H^n}f(\mathcal C(z,t)){\rm Jac}_{\mathcal C}(z,t){\rm d}z{\rm d}t.
\end{equation}
 If $w=(z,t)$, $v=(z',t')$ and $\zeta=\mathcal C(w)$, $\eta=\mathcal C(v)$, one has that 
 \begin{equation}\label{tavolsagok}
 d_S(\zeta,\eta)=d_{KC}(w,v)\left(\frac{{4}}{((1+|z|^2)^2+t^2)}\right)^\frac{1}{4}\left(\frac{{4}}{((1+|z'|^2)^2+(t')^2)}\right)^\frac{1}{4}.
 \end{equation}
 
 %	\subsection{CR fractional sub-Laplacian operator on $\mathbb H^n$} 

 \subsection{Spherical/zonal harmonics on $S^{2n+1}$}  The Hilbert space $L^2(S^{2n+1})$, endowed by the inner product 
 $$(U,V)=\int_{S^{2n+1}}U\overline{V}{\rm d}\eta,$$ can be decomposed into ${\bf U}(n+1)$-irreducible components as 
 $$L^2(S^{2n+1})= \bigoplus_{j,k\geq 0}\mathcal H_{j,k},$$ where $\mathcal H_{j,k}$ denotes the space of harmonic polynomials $p(z,\overline z)$ on $\mathbb C^{n+1}$ restricted to $S^{2n+1}$ that are homogeneous of degree $j$ and $k$ in the variables $z$ and $\overline z$, respectively. We notice that the dimension of $\mathcal H_{j,k}$ is 
 $$m_{j,k}=\frac{(j+n-1)!(k+n-1)!(j+k+n)}{n!(n-1)!j!k!}.$$
 Moreover, if $\{Y_{j,k}^l\}_{l=\overline{1,m_{j,k}}}$ is an orthonormal  basis of $\mathcal H_{j,k},$ then the zonal harmonics are defined by 
 %\begin{equation}\label{zonal-harmonics}
\begin{equation}\label{zonal-0}
 \Phi_{j,k}(\zeta,\eta)=\sum_{l=1}^{m_{j,k}}Y_{j,k}^l(\zeta)\overline{Y_{j,k}^l(\eta)},\ \ \zeta,\eta\in S^{2n+1}.
  \end{equation}
 %\end{equation}
 We recall that $ \Phi_{j,k}$ can be represented as
 \begin{equation}\label{zonal}
  \Phi_{j,k}(\zeta,\eta)=\frac{(\max\{j,k\}+n-1)!(j+k+n)}{\omega_{2n+1}n!(\max\{j,k\})!}\langle\zeta,\overline \eta\rangle^{|j-k|}P_k^{(n-1,|j-k|)}(2\langle\zeta,\overline \eta\rangle^2-1),
 \end{equation}
 where $P_k^{(n,l)}$ denotes the Jacobi polynomials and $\omega_{2n+1}$ is the usual measure of $S^{2n+1}$.  
 
 \subsection{Fractional Sobolev spaces on $S^{2n+1}$ and $\mathbb H^n$}\label{subsection2-2} 
 The usual sub-Laplacian  on $\mathbb H^n$ is defined as 
 $$\mathcal L=-\frac{1}{4}\sum_{j=1}^n(X_j^2+Y_j^2).$$
If we introduce the differential operators $$T_j=\frac{\partial}{\partial \eta_j}-\overline \eta_j\sum_{k=1}^{n+1}\eta_k\frac{\partial}{\partial \eta_k}, \ \ \ \ \overline {T_j}=\frac{\partial}{\partial \overline{\eta_j}}- \eta_j\sum_{k=1}^{n+1}\overline{\eta_k}\frac{\partial}{\partial \overline{\eta_k}},\ \ j\in\{1,...,n+1\},$$
%and $$\mathcal L_0=-\frac{1}{2}\sum_{j=1}^{n+1}(T_j\overline {T_j}+\overline {T_j}T_j),$$
 the  conformal sub-Laplacian on $S^{2n+1}$ is given by $$\mathcal D=-\frac{1}{2}\sum_{j=1}^{n+1}(T_j\overline {T_j}+\overline {T_j}T_j)+\frac{n^2}{4}.$$ 
   Note that for every $Y_{j,k}\in \mathcal H_{j,k}$ one has 
   $$\mathcal DY_{j,k}=\lambda_j\lambda_k Y_{j,k},$$ 
   where $\lambda_j=j+\frac{n}{2}.$ 
   
  Let $0<\gamma<\frac{Q}{2}=n+1$ be fixed. Given $U\in L^2(S^{2n+1})$, its Fourier representation is  
   $$U=\sum_{j,k\geq  0}\sum_{l=1}^{m_{j,k}}c_{j,k}^l(U)Y_{j,k}^l$$ with Fourier coefficients 
   $c_{j,k}^l(U)=\displaystyle\int_{S^{2n+1}}UY_{j,k}^l{\rm d}\eta$. Accordingly, we may define 
   $$\mathcal D^\frac{\gamma}{2}U=\sum_{j,k\geq  0}\sum_{l=1}^{m_{j,k}}(\lambda_j\lambda_k)^\frac{\gamma}{2}c_{j,k}^l(U)Y_{j,k}^l.$$
   The fractional Sobolev space over $S^{2n+1}$ is defined as
   $$H^\gamma(S^{2n+1})=W^{\gamma,2}(S^{2n+1})=\left\{U\in L^2(S^{2n+1}):\mathcal D^\frac{\gamma}{2}U\in L^2(S^{2n+1}) \right\},$$
   endowed with the inner product and norm 
   $$(U,V)_{H^\gamma}=\int_{S^{2n+1}}\mathcal D^\frac{\gamma}{2}U\overline{\mathcal D^\frac{\gamma}{2} V}{\rm d}\eta\ \ {\rm and}\ \ \|U\|_{H^\gamma}=(U,U)_{H^\gamma}^\frac{1}{2}= \left(\sum_{j,k\geq  0}\sum_{l=1}^{m_{j,k}}(\lambda_j\lambda_k)^{\gamma}|c_{j,k}^l(U)|^2\right)^\frac{1}{2}.$$ The  norm $\|\cdot\|_{H^\gamma}$ is equivalent to the norm coming from the inner product $$(U,V)_\gamma=\sum_{j,k\geq  0}\sum_{l=1}^{m_{j,k}}\lambda_j(\gamma)\lambda_k(\gamma)c_{j,k}^l(U)\overline{c_{j,k}^l(V)},$$
   where $$\lambda_j(\gamma)=\frac{\Gamma\left(\frac{Q+2\gamma}{4}+j\right)}{\Gamma\left(\frac{Q-2\gamma}{4}+j\right)},\ \ j\in \mathbb N_0=\{0,1,2,...\};$$  indeed, by asymptotic approximation of the Gamma function $\Gamma$, one has $\lambda_j(\gamma)\sim j^\gamma$. 

The intertwining operator $\mathcal A_\gamma$ of order $2\gamma$ on $S^{2n+1}$ is given by 
$${\rm Jac}_\tau^\frac{Q+2\gamma}{2Q}(\mathcal A_\gamma U)\circ \tau=\mathcal A_\gamma({\rm Jac}_\tau^\frac{Q-2\gamma}{2Q}(U\circ \tau))\ \ {\rm for\ all}\ \tau\in {\rm Aut}(S^{2n+1}),U\in  C^\infty(S^{2n+1}),$$ 
where ${\rm Aut}(S^{2n+1})$ and ${\rm Jac}_\tau$ denote the group of automorphisms on $S^{2n+1}$ and the Jacobian of $\tau\in {\rm Aut}(S^{2n+1})$, respectively.
%16 oldal -BFM 
In fact, the latter definition can be extended to every $U\in H^\gamma(S^{2n+1})$. 
Note that $\mathcal A_\gamma$ may by characterized (up to a constant) by its action on $\mathcal H_{j,k}$ as 
\begin{equation}\label{sajat-ertek-A}
\mathcal A_\gamma Y_{j,k}={\lambda_j(\gamma)\lambda_k(\gamma)}%{c_\gamma}
Y_{j,k}, \ \ Y_{j,k}\in \mathcal H_{j,k}.
\end{equation}
%where $c_\gamma=???$ 
Therefore, 
\begin{equation}\label{osszefugges}
(U,V)_{\gamma}=\int_{S^{2n+1}}\overline  V \mathcal A_\gamma U  {\rm d}\eta.
\end{equation}  
In particular, $\lambda_j(1)=\lambda_j$ for every $j\in \mathbb N_0$ and  
$\mathcal A_1=\mathcal D$. Moreover, according to Frank and Lieb \cite{FL}, for every real-valued function $U\in H^\gamma(S^{2n+1})$, one has the sharp fractional Sobolev inequality on the CR sphere $S^{2n+1}$, i.e.,
\begin{equation}\label{sharp-FL}
\left(\int_{S^{2n+1}}|U(\eta)|^\frac{2Q}{Q-2\gamma}{\rm d}\eta\right)^\frac{Q-2\gamma}{Q}\leq C(\gamma,n)\int_{S^{2n+1}}U(\eta)\mathcal A_\gamma U(\eta){\rm d}\eta,
\end{equation}
where $$C(\gamma,n)=\frac{\Gamma\left(\frac{n+1-\gamma}{2}\right)^2}{\Gamma\left(\frac{n+1+\gamma}{2}\right)^2}\omega_{2n+1}^{-\frac{\gamma}{n+1}}.$$
%and all the optimizers are translates, dilates or constant multiples of the function 
 
  The 	CR fractional sub-Laplacian operator on $\mathbb H^n$ is defined by 
 $$\mathcal L_\gamma=|2T|^\gamma\frac{\Gamma\left(\mathcal L|2T|^{-1}+\frac{1+\gamma}{2}\right)}{\Gamma\left(\mathcal L|2T|^{-1}+\frac{1-\gamma}{2}\right)}.$$
 Direct computation shows that
 $\mathcal L_1=\mathcal L$, $\mathcal L_2=\mathcal L^2- |T|^2$. 
%  and 
% \begin{equation}
% \mathcal L_\gamma(su)=s^\gamma\mathcal L_\gamma u,\ \ \forall s\geq 0, u\in C^\infty(\mathbb H^n).
% \end{equation}
  Moreover, the relationship between $\mathcal L_\gamma$ and $\mathcal A_\gamma$ is given by 
 \begin{equation}\label{A-es-L}
 \mathcal L_\gamma\left((2{\rm Jac}_{\mathcal C})^\frac{Q-2\gamma}{2Q}(U\circ \mathcal C)\right)=(2{\rm Jac}_{\mathcal C})^\frac{Q+2\gamma}{2Q}(\mathcal A_\gamma U)\circ \mathcal C,\ \ \ \forall U\in H^\gamma(S^{2n+1}). 
 \end{equation}
 % 21 oldal BFM
% $${\rm Jac}_h^\frac{Q+2\gamma}{2Q}(\mathcal L_\gamma v)\circ h=\mathcal L_\gamma({\rm Jac}_h^\frac{Q-2\gamma}{2Q}(v\circ h))\ \ {\rm for\ all}\ h\in {\rm Aut}(\mathbb H^n),v\in  C^\infty(\mathbb H^n),$$ 
% where ${\rm Aut}(S^{2n+1})$ and ${\rm Jac}_\tau$ denote the group of automorphisms on $S^{2n+1}$ and the Jacobian of $\tau\in {\rm Aut}(S^{2n+1})$, respectively.

  The fractional Sobolev space over $\mathbb H^n$ % (for real-valued functions) 
  is defined by  
 $$D^{\gamma}({\mathbb{H}^n})=\left\{u\in L^{\frac{2Q}{Q-2\gamma}}(\mathbb H^n): a_\gamma[u]<+\infty\right\},$$
 where the quadratic form
  $a_\gamma$ is associated to the operator $\mathcal L_\gamma$, i.e.,
  $$a_\gamma[u]=\int_{\mathbb H^n} \overline u \mathcal L_\gamma u {\rm d}z{\rm d}t.$$ 
  The  form $a_\gamma$ can be equivalently represented by means of spectral decomposition, see \cite[p. 126]{Frank-Adv-Math}.

% The inner product on $H^\gamma(S^{2n+1})$ is given by 
% \begin{equation}\label{inner-sphere}
%(U,V)_\#=2^{\gamma-\frac{Q}{2}}\sum_{j,k\geq  0}\sum_{l=0}^{m_{jk}}\lambda_j(\gamma)\lambda_k(\gamma)c_{jk}^l(U)\overline{c_{jk}^l(V)}.
% \end{equation}

%\textbf{	For simplicity, we denote by $2^*_\gamma=\frac{2Q}{Q-2\gamma}$ the fractional critical exponent, and 
%	}

	\section{Proof of Main Theorem}\label{section-3}

	\subsection{Equivalent critical problems on $\mathbb H^n$ and  $S^{2n+1}$.}\label{section-3-1} Let $\gamma\in (0,n+1)$  be fixed. We consider the fractional Yamabe problem on the CR sphere, i.e., 
		\begin{equation*}
	\left\{
	\begin{array}{lll}
	\mathcal A_\gamma U=|U|^\frac{4\gamma}{Q-2\gamma}U & \mbox{on} & S^{2n+1}, \\
	U\in H^\gamma(S^{2n+1}).
	\end{array}%
	\right.\eqno{{(\textbf{FYS})_\gamma}}
	\end{equation*}
	%where the operator $\mathcal A_\gamma$ and Sobolev space $H^\gamma(S^{2n+1})$ are introduced in \S \ref{subsection2-2}. 
	Hereafter, we are considering real-valued functions both in $H^\gamma(S^{2n+1})$ and $D^\gamma(\mathbb H^n)$, respectively.
	The main result of this subsection constitutes the bridge between {\rm ${(\textbf{FYS})_\gamma}$} and {\rm ${(\textbf{FYH})_\gamma}$} as follows:  
	
	\begin{proposition}\label{weak-equivalent}
	 Let $0<\gamma<\frac{Q}{2}=n+1$. Then	$U\in H^\gamma(S^{2n+1})$ is a weak solution of {\rm ${(\textbf{FYS})_\gamma}$} if and only if 
%		$(z,t)\mapsto u(z,t)=2^\alpha {\rm Jac}_{\mathcal C}(z,t)^\frac{Q-2\gamma}{2Q}U(\mathcal C(z,t))$ 
$u=(2{\rm Jac}_{\mathcal C})^\frac{Q-2\gamma}{2Q}U\circ \mathcal C\in  D^\gamma(\mathbb H^n)$
		is a weak solution of {\rm${(\textbf{FYH})_\gamma}$}.  
	\end{proposition}

{\it Proof.} We first prove the following 

\noindent \textit{\underline{Claim}}: {\it Let $U:S^{2n+1}\to \mathbb R$ and $u:\mathbb H^n\to \mathbb R$ be two functions such that $u={\rm Jac}_{\mathcal C}^\frac{Q-2\gamma}{2Q}U\circ \mathcal C$. Then $U\in H^\gamma(S^{2n+1})$ if and only if $u\in  D^\gamma(\mathbb H^n)$.}

Fix $U\in H^\gamma(S^{2n+1})$; we shall prove first that  $(z,t)\mapsto u(z,t)= {\rm Jac}_{\mathcal C}(z,t)^\frac{Q-2\gamma}{2Q}U(\mathcal C(z,t))$ belongs to $ D^\gamma(\mathbb H^n)$. By (\ref{valtozo-csere}) one has 
\begin{eqnarray}\label{crit-estim}
\nonumber \int_{\mathbb H^n}|u(z,t)|^\frac{2Q}{Q-2\gamma}{\rm d}z{\rm d}t&=&\int_{\mathbb H^n}{\rm Jac}_{\mathcal C}(z,t)|U(\mathcal C(z,t))|^\frac{2Q}{Q-2\gamma}{\rm d}z{\rm d}t\\&=&\int_{S^{2n+1}}|U(\eta)|^\frac{2Q}{Q-2\gamma}{\rm d}\eta.
\end{eqnarray}
Furthermore, by the fractional Sobolev inequality (\ref{sharp-FL}) and relation (\ref{sajat-ertek-A}), one has that 
\begin{eqnarray*}
\left(\int_{S^{2n+1}}|U(\eta)|^\frac{2Q}{Q-2\gamma}{\rm d}\eta\right)^\frac{Q-2\gamma}{Q}&\leq& C(\gamma,n)\int_{S^{2n+1}}U(\eta)\mathcal A_\gamma U(\eta){\rm d}\eta\\&=&C(\gamma,n)\sum_{j,k\geq  0}\sum_{l=1}^{m_{j,k}}\lambda_j(\gamma)\lambda_k(\gamma)|c_{j,k}^l(U)|^2\\&\leq &C'(\gamma,n)\|U\|_{H^\gamma}^2\\&<&+\infty,
\end{eqnarray*}
where $C'(\gamma,n)=C_\gamma C(\gamma,n)$ and $C_\gamma>0$ is such that $(V,V)_\gamma\leq C_\gamma\|V\|_{H^\gamma}^2$ for every $V\in H^\gamma(S^{2n+1})$; thus $u\in  L^{\frac{2Q}{Q-2\gamma}}(\mathbb H^n)$. Moreover, by (\ref{A-es-L}) and (\ref{valtozo-csere}) one has
\begin{eqnarray}\label{ez-is-kell}
\nonumber a_\gamma[u]&=&\int_{\mathbb H^n}u \mathcal L_\gamma u {\rm d}z{\rm d}t= \nonumber \\&=& \nonumber
2^{\alpha'} \int_{\mathbb H^n} {\rm Jac}_{\mathcal C}(z,t)^\frac{Q-2\gamma}{2Q}U(\mathcal C(z,t)) \mathcal L_\gamma \left((2{\rm Jac}_{\mathcal C}(z,t))^\frac{Q-2\gamma}{2Q}U(\mathcal C(z,t))\right){\rm d}z{\rm d}t\\
\nonumber &=&2^{\alpha'} \int_{\mathbb H^n} {\rm Jac}_{\mathcal C}(z,t)^\frac{Q-2\gamma}{2Q}U(\mathcal C(z,t))  (2{\rm Jac}_{\mathcal C}(z,t))^\frac{Q+2\gamma}{2Q}(\mathcal A_\gamma U)(\mathcal C(z,t)){\rm d}z{\rm d}t \nonumber
\\
&=&2^{\alpha''} \int_{\mathbb H^n} U(\mathcal C(z,t))(\mathcal A_\gamma U)(\mathcal C(z,t))  {\rm Jac}_{\mathcal C}(z,t){\rm d}z{\rm d}t \nonumber
\\
&=&2^{\alpha''} \int_{S^{2n+1}} U(\eta)\mathcal A_\gamma U(\eta)  {\rm d}\eta\\
&<&+\infty,\nonumber
\end{eqnarray}
where $\alpha'=-\frac{Q-2\gamma}{2Q}$ and $\alpha''=\alpha'+\frac{Q+2\gamma}{2Q}=\frac{2\gamma}{Q}.$ Therefore, $u\in D^\gamma(\mathbb H^n)$. 

Conversely, let us assume that $u\in  D^\gamma(\mathbb H^n)$. In particular, we have that $u\in L^{\frac{2Q}{Q-2\gamma}}(\mathbb H^n)$, thus by relation (\ref{crit-estim}) it turns out that $U\in L^{\frac{2Q}{Q-2\gamma}}(S^{2n+1})$; therefore, $U\in L^2(S^{2n+1})$. Furthermore, by (\ref{ez-is-kell}) we also have that
$$ \int_{S^{2n+1}} U(\eta)\mathcal A_\gamma U(\eta)  {\rm d}\eta=2^{-\alpha''}a_\gamma[u]<+\infty,$$
i.e.,  $U\in H^\gamma(S^{2n+1})$, which concludes the proof of \textit{Claim}.

Let $U\in H^\gamma(S^{2n+1})$ be a weak solution of {\rm ${(\textbf{FYS})_\gamma}$}; then we have 
\begin{equation}\label{rweak-1}
\int_{S^{2n+1}}\mathcal A_\gamma U V{\rm d}\eta=\int_{S^{2n+1}}|U|^\frac{4\gamma}{Q-2\gamma}UV{\rm d}\eta\ \ \ {\rm for\ every}\ V\in H^\gamma(S^{2n+1}).
\end{equation}
Let $v\in  D^\gamma(\mathbb H^n)$ be  arbitrarily fixed and define $V=({\rm Jac}_{\mathcal C}\circ \mathcal C^{-1})^\frac{2\gamma-Q}{2Q}v\circ \mathcal C^{-1}.$ Since $v= {\rm Jac}_{\mathcal C}^\frac{Q-2\gamma}{2Q}V\circ \mathcal C$, by the \textit{Claim} we have that $V\in H^\gamma(S^{2n+1})$. Accordingly, the function $V$ can be used as a test-function in (\ref{rweak-1}), obtaining 
$$ \int_{S^{2n+1}}\mathcal A_\gamma U ({\rm Jac}_{\mathcal C}\circ \mathcal C^{-1})^\frac{2\gamma-Q}{2Q}v\circ \mathcal C^{-1}{\rm d}\eta=\int_{S^{2n+1}}|U|^\frac{4\gamma}{Q-2\gamma}U({\rm Jac}_{\mathcal C}\circ \mathcal C^{-1})^\frac{2\gamma-Q}{2Q}v\circ \mathcal C^{-1}{\rm d}\eta.$$
By a change of variables, it follows that 
$$ \int_{\mathbb H^n}(\mathcal A_\gamma U\circ \mathcal C ) ({\rm Jac}_{\mathcal C})^{\frac{2\gamma-Q}{2Q}+1}v{\rm d}z{\rm d}t =\int_{\mathbb H^n}|U\circ \mathcal C|^\frac{4\gamma}{Q-2\gamma}(U\circ \mathcal C)({\rm Jac}_{\mathcal C})^{\frac{2\gamma-Q}{2Q}+1}v{\rm d}z{\rm d}t.$$
This relation and (\ref{A-es-L}) imply that 
$$ 2^{-\frac{Q+2\gamma}{2Q}}\int_{\mathbb H^n}\mathcal L_\gamma\left((2{\rm Jac}_{\mathcal C})^\frac{Q-2\gamma}{2Q}(U\circ \mathcal C)\right) v{\rm d}z{\rm d}t =\int_{\mathbb H^n}|U\circ \mathcal C|^\frac{4\gamma}{Q-2\gamma}U\circ \mathcal C({\rm Jac}_{\mathcal C})^{\frac{2\gamma+Q}{2Q}}v{\rm d}z{\rm d}t.$$
Since $u= (2{\rm Jac}_{\mathcal C})^\frac{Q-2\gamma}{2Q}U\circ \mathcal C$, the latter equality is equivalent to 
$$ \int_{\mathbb H^n}\mathcal L_\gamma u v{\rm d}z{\rm d}t =\int_{\mathbb H^n}|u|^\frac{4\gamma}{Q-2\gamma}uv{\rm d}z{\rm d}t,$$
which means precisely that $u\in  D^\gamma(\mathbb H^n)$ is a weak solution of {\rm ${(\textbf{FYH})_\gamma}$}. The converse argument works in a similar way. \hfill $\square$

%By relations (\ref{sajat-ertek-A}), it turns out that 
%\begin{eqnarray*}
%\int_{S^{2n+1}}\mathcal A_\gamma W V{\rm d}\eta&=&\sum_{j,k\geq  0}\sum_{l=0}^{m_{jk}}c_{jk}^l(W)\overline{c_{jk}^l(V)} \int_{S^{2n+1}}\mathcal A_\gamma Y_{jk}^l \overline{Y_{jk}^l}{\rm d}\eta\\&=&\sum_{j,k\geq  0}\sum_{l=0}^{m_{jk}}c_{jk}^l(W)\overline{c_{jk}^l(V)} \int_{S^{2n+1}}\mathcal A_\gamma Y_{jk}^l \overline{Y_{jk}^l}{\rm d}\eta
%\end{eqnarray*}

\begin{remark}\rm \label{weak-2}  One can provide an alternative proof to Proposition \ref{weak-equivalent} by exploring the explicit form of the fundamental solution of $\mathcal L_\gamma$; a similar approach is due to  Bartsch, Schneider and Weth \cite{BSW} for the polyharmonic operator $(-\Delta)^m$ in $\mathbb R^N$, where $m\in \mathbb N$ and $N>2m$. For completeness, we sketch the proof. 
	
We recall that  the fundamental solution of $\mathcal L_\gamma$ is  
\begin{equation}\label{fundamental}
\mathcal  L_\gamma^{-1}((z,t),(z',t'))=\frac{c_\gamma}{2}d_{KC}^{2\gamma-Q}((z,t),(z',t')),
\end{equation}
where $$c_\gamma=\frac{2^{n-\gamma}\Gamma\left(\frac{Q-2\gamma}{4}\right)^2}{\pi^{n+1}\Gamma(\gamma)},$$ see Branson, Fontana and Morpurgo \cite[p. 21]{BFM}. 
For every  $\psi\in L^{\frac{2Q}{Q+2\gamma}}(S^{2n+1})$ we introduce  the function 
		\begin{equation}\label{K-definictio}
		[\mathcal K_\gamma \psi](\zeta)=c_\gamma \int_{S^{2n+1}}\psi(\eta)|1-\langle \zeta, \overline \eta \rangle|^\frac{2\gamma-Q}{2}{\rm d}\eta.
		\end{equation}
	One can prove that $\mathcal K_\gamma \psi\in H^\gamma(S^{2n+1})$ for every $\psi\in L^{\frac{2Q}{Q+2\gamma}}(S^{2n+1})$. Moreover, the Funk-Hecke theorem on the CR sphere $S^{2n+1}$ gives
		$$[\mathcal K_\gamma Y_{j,k}](\zeta)=\frac{2^{\frac{Q}{2}-\gamma}}{{\lambda_j(\gamma)\lambda_k(\gamma)}}Y_{j,k}(\zeta),$$
		see Frank and Lieb \cite[Corollary 5.3]{FL}. 
		Thus, a direct computation yields that 
		$$(\mathcal K_\gamma \psi,V)_{\gamma}=2^{\frac{Q}{2}-\gamma}\int_{S^{2n+1}} \psi  V  {\rm d}\eta\ \ \ {\rm for\ every}\ V\in H^\gamma(S^{2n+1}).$$
		 Note that if $U\in H^\gamma(S^{2n+1})$ is a weak solution of {\rm ${(\textbf{FYS})_\gamma}$}, the latter relation implies that 
		\begin{equation}\label{K-atirva}
		\mathcal K_\gamma(|U|^\frac{4\gamma}{Q-2\gamma}U)=2^{\frac{Q}{2}-\gamma}U\ \ {\rm on}\ \ S^{2n+1}.
		\end{equation}
		 Accordingly,  by  relations (\ref{K-atirva}), (\ref{K-definictio}), (\ref{valtozo-csere}) and  (\ref{tavolsagok}), it turns out  that
		\begin{eqnarray*}
		u(z,t)&=&(2{\rm Jac}_{\mathcal C}(z,t))^\frac{Q-2\gamma}{2Q}U(\mathcal C(z,t))\\&=&2^{-\frac{Q}{2}+\gamma}(2{\rm Jac}_{\mathcal C}(z,t))^\frac{Q-2\gamma}{2Q}\mathcal K_\gamma(|U(\mathcal C(z,t))|^\frac{4\gamma}{Q-2\gamma}U(\mathcal C(z,t)))\\
		&=&\frac{c_\gamma}{2} \int_{\mathbb H^n} d_{KC}^{2\gamma-Q}((z,t),(z',t')|u(z',t')|^\frac{4\gamma}{Q-\gamma}u(z',t'){\rm d}z'{\rm d}t',\ \ \ (z,t)\in \mathbb H^n.
		\end{eqnarray*}	
The latter relation is equivalent to the fact that
\begin{equation}\label{convol}
u(z,t)=\frac{c_\gamma}{2}(|u|^\frac{4\gamma}{Q-\gamma}u)*d_{KC}^{2\gamma-Q}((z,t),\cdot),\ \ \ (z,t)\in \mathbb H^n,
\end{equation}
 where $'*'$ denotes the (noncommutative) convolution operation on the Heisenberg group $\mathbb H^n$. By (\ref{fundamental}), a similarly argument as in Folland \cite[Theorem 2]{Folland-BAMS} gives that $\mathcal L_\gamma u=|u|^\frac{4\gamma}{Q-\gamma}u$ on $\mathbb H^n,$	which concludes the claim. 	
\end{remark}

%-----
%
%The proof is divided into \textbf{three} steps. 
%
%\textit{Step 1.}  
%
%
%\begin{equation*}
%\left\{
%\begin{array}{lll}
%\mathcal L_\gamma w=|w|^\frac{4\gamma}{Q-2\gamma}w & \mbox{on} & \mathbb H^n, \\
%w\in \mathcal D^{\gamma}({\mathbb{H}^n}).
%\end{array}%
%\right.\eqno{{(\textbf{FYH})_\gamma}}
%\end{equation*}
%
%
%$$D^{\gamma}({\mathbb{H}^n})=\{w\in L^{2^*_\gamma}(\mathbb H^n): a_\gamma(w)<+\infty\},$$
%where 
%$2^*_\gamma=\frac{2Q}{Q-2\gamma}$ and $a_\gamma$ is the quadratic form associated to the operator $\mathcal L_\gamma$, see ???. 

	\medskip
	
	\subsection{Compactness}\label{section-3-2} According to Frank and Lieb \cite{FL}, see also (\ref{sharp-FL}), the embedding $H^\gamma(S^{2n+1})\hookrightarrow L^\frac{2Q}{Q-2\gamma}(S^{2n+1})$ is continuous, but not compact. This subsection is devoted to regain certain compactness by using suitable group actions on the CR sphere $S^{2n+1}$.  
	
	To complete this purpose, let $n_j\in \mathbb N$, $j=1,...,k,$  with $n_1+...+n_k=n+1.$ 
Associated to these numbers, let  
\begin{equation}\label{splitting}
G={\bf U}(n_1)\times...\times {\bf U}(n_k)
\end{equation}
	be the subgroup of the unitary group ${\bf U}(n+1)=\left\{g\in {\bf O}(2n+2):gJ=Jg\right\},$ where $$J=\left[\begin{matrix}
	0 & I_{\mathbb R^{n+1}}\\
	-I_{\mathbb R^{n+1}}& 0
	\end{matrix}\right].$$ Let
	$$H_G^\gamma(S^{2n+1})=\left\{U\in H^\gamma(S^{2n+1}): g\circ U=U\ \ \ {\rm for\ every}\ g\in G \right\}$$
	 be the subspace of $G$-invariant functions of $H^\gamma(S^{2n+1})$, where
	 \begin{equation}\label{elso-action}
	 (g\circ U)(\eta)=U(g^{-1}\eta),\ \  \eta\in S^{2n+1}.
	 \end{equation}
	 It is clear that $H_G^\gamma(S^{2n+1})$   is  an 
	 infinite-dimensional	 closed	 subspace	 of $H^\gamma(S^{2n+1})$, whenever $k\geq 2$ in the  splitting (\ref{splitting}). 
	
	\medskip
	With the above notations in our mind, a Ding-Hebey-Vaugon-type compactness result  reads as follows:  
	\begin{proposition}\label{proposition-compact}
		Let $\gamma\in \bigcup_{k=1}^n[k,\frac{kQ}{Q-1})$. The embedding $H_G^\gamma(S^{2n+1})\hookrightarrow L^\frac{2Q}{Q-2\gamma}(S^{2n+1})$ is compact, where $G={\bf U}(n_1)\times...\times {\bf U}(n_k)$ is any choice with $n_j\in \mathbb N$, $j=1,...,k,$  and $n_1+...+n_k=n+1.$ 
	\end{proposition}
	 
	 {\it Proof.} First, when $G={\bf U}(n+1)$, the space $H_G^\gamma(S^{2n+1})$ contains precisely the constant functions defined on $S^{2n+1};$ in this case the proof is trivial. 
%	 We assume that $n_1\leq n_2$, thus $0<\gamma<n_2.$ Let $\tilde Q=2n_2$; we are going to prove that the embedding $H_G^\gamma(S^{2n+1})\hookrightarrow L^p(S^{2n+1})$ is continuous for $1\leq p\leq \frac{2\tilde Q}{\tilde Q-2\gamma}$, and compact for $1\leq p< \frac{2\tilde Q}{\tilde Q-2\gamma}$. 
%	   
%	 In order to prove this, let us fix $U\in H_G^\gamma(S^{2n+1})$. In particular, for $g=(g_1,g_2)\in G$ with $g_i\in {\bf U}(n_i)$, the action of $G$ on $S^{2n+1}$ is given by $g\eta=g(z_1,z_2)=(g_1z_1,g_2z_2)$, where $z_i\in \mathbb C^{n_i}$, $i\in \{1,2\}$ and $\eta=(z_1,z_2)\in S^{2n+1}$.  Since $U\in H_G^\gamma(S^{2n+1})$, it turns out that $U$ depends on $|\eta_2|$ (or equivalently, on $|\eta_1|$). Therefore, there exists a bijective correspondence $\eta=(z_1,z_2)\mapsto (\eta_1,\eta_2,\theta)\in S^{2n_1-1}\times S^{2n_2-1}\times [0,2\pi)$ where $\eta=\sin \theta \eta_1+\cos \theta \eta_2$; in particular, $\eta_i=\frac{z_i}{|z_i|},$ $i\in \{1,2\}$ (where the trivial cases $z_1=0$ or $z_2=0$ are treated in the usual way). 
%	 
%	  $U(\eta)=U(|\eta_1|,|\eta_2|)$, which means that . Let $\tilde U(\eta_2):=U(|\eta_1|,|\eta_2|)=U(\eta).$  
%	
%	Since $\frac{2Q}{Q-2\gamma}< \frac{2\tilde Q}{\tilde Q-2\gamma}$, our claim will follow at once from the above statement.
 
 In the general case, we recall by Maalaoui and  Martino \cite[Lemma 3.3]{Martino-0} that the embedding  $W^{1,2}_G(S^{2n+1}) = H_G^1(S^{2n+1})\hookrightarrow L^{q}(S^{2n+1})$ is compact for every $1\leq q<q_1^*$, where $q_1^*=\frac{2(Q-1)}{Q-3}$ is the Riemannian critical exponent on the $(Q-1)-$dimensional sphere $S^{2n+1}$. 
 
 By our assumption
 $\gamma\in \bigcup_{k=1}^n[k,\frac{kQ}{Q-1})$
 % $0<\gamma<\frac{Q}{2}$ and $[\gamma(1-\frac{1}{Q})]+1=[\gamma]$ we have that 
 %$0<\gamma-\gamma(1-\frac{1}{Q})=\frac{\gamma}{Q}<\frac{1}{2}$, 
 we have that $l:=[\gamma]\geq 1$ and 
\begin{equation}\label{gamma-ra-felt}
\gamma\left(1-\frac{1}{Q}\right)< l\leq \gamma.
\end{equation}
The iterative argument developed by Aubin \cite[Proposition 2.11]{Aubin}, applied for $l$ times,  gives that the embedding $W^{l,2}_G(S^{2n+1}) = H_G^l(S^{2n+1})\hookrightarrow L^{q}(S^{2n+1})$ is compact for every $1\leq q<q_l^*$, where $q_l^*=\frac{2(Q-1)}{Q-1-2l}$. On one hand, since $l\leq \gamma$, we have that $H^\gamma_G(S^{2n+1})=W^{\gamma,2}_G(S^{2n+1})\subset W^{l,2}_G(S^{2n+1})$. 
On the other hand, the left hand side of (\ref{gamma-ra-felt}) is equivalent to $q_l^*> \frac{2Q}{Q-2\gamma}$. Combining these facts, we have the chain of inclusions 
$$H^\gamma_G(S^{2n+1})\subset W^{l,2}_G(S^{2n+1})\hookrightarrow L^\frac{2Q}{Q-2\gamma}(S^{2n+1}),$$ 
where the latter embedding is compact. 
	 \hfill $\square$
	 
	 \begin{remark}\rm\label{remark-kell}
	 	Our  assumption $\gamma\in \bigcup_{k=1}^n[k,\frac{kQ}{Q-1})$ is important to guarantee the left hand side of (\ref{gamma-ra-felt}), which in turn, implies that $\frac{2Q}{Q-2\gamma}$ is within the range $[1,q_l^*)$ where the embedding $W^{l,2}_G(S^{2n+1}) \hookrightarrow L^{q}(S^{2n+1})$ is compact, $q\in [1,q_l^*)$. We are wondering if this assumption can be removed  in order to prove the compactness of the above embedding for the whole spectrum $(0,\frac{Q}{2})$ of the parameter $\gamma.$
	 \end{remark}

\medskip

	\subsection{Special group actions}\label{section-3-3} 
 The goal of this subsection is to describe
	symmetrically different functions belonging to
	$H^\gamma(S^{2n+1})$ via subgroups of the form $G={\bf U}(n_1)\times...\times {\bf U}(n_k)$  with $n_1+...+n_k=n+1$. 
	 To handle this problem,  we explore a
	Rubik-cube technique, described in a slightly different manner in Balogh and Krist\'aly \cite{BK}; roughly speaking,   $n+1$
	corresponds to the number of total sides of the cube, while the sides of
	the cube are certain blocks in the decomposition subgroup
	$G={\bf U}(n_1)\times...\times {\bf U}(n_k)$.

	To be more precise, let $n\geq 1$ and for $i\in \{1,...,\left[\frac{n+1}{2}\right]\}$ we
	consider the subgroup of the unitary group ${\bf U}(n+1)$ as
	$$G_{i}=\left\{
	\begin{array}{lcl}
	\left[
	\begin{array}{cc}
	{\bf U}(\frac{n+1}{2}) & 0  \\
	0 & {\bf U}(\frac{n+1}{2})  \\
	\end{array}
	\right],& \mbox{if}&
	n+1=2i,\\
	\\
\left[
\begin{array}{ccc}
{\bf U}(i) & 0 & 0  \\
0& {\bf U}(n+1-2i) & 0 \\
0& 0 & {\bf U}(i)
\end{array}
\right],& \mbox{if}& n+1\neq 2i.
	\end{array}\right.$$
	It is clear that a particular $G_i$ does not act transitively on the sphere $S^{2n+1}$. However, to recover the transitivity, we shall combine different groups of the type $G_i$; for further use, let $ [ G_i;G_j]$ be the group
	generated by $G_i$ and $G_j.$
	
	\begin{lemma}\label{Rubik}
		Let $i,j\in \{1,...,\left[\frac{n+1}{2}\right]\}$ with $i\neq j.$ Then
		the group $[G_i;G_j]$ acts transitively on the CR sphere
		${S}^{2n+1}$. 
	\end{lemma}
	
	{\it Proof.}  Without loss of generality, we assume that $j>i$. For further use, let $0_k=%(\underbrace{0,...,0}_{k\ {\rm times}}).
	(0,...,0)\in \mathbb C^k=\mathbb R^{2k}$, $k\in \{1,...,n\}.$ Let us fix $\eta=(\eta_1,\eta_2,\eta_3)\in {S}^{2n+1}$
	arbitrarily with $\eta_1,\eta_3\in \mathbb C^j$ and $\eta_2\in \mathbb
	C^{n+1-2j}$; clearly,  $\eta_2$  disappears from $\eta$ whenever $2j=n+1$.
	Taking into account the fact that ${\bf U}(j)$ acts transitively on $S^{2j-1}$, there are
	$g_j^1,g_j^2\in {\bf U}(j)$ such that if $g_j= g_j^1 \times
	I_{\mathbb C^{n+1-2j}}\times g_j^2\in G_j,$ then
	$g_j\eta=(0_{j-1},0,|\eta_1|,\eta_2,|\eta_3|,0,0_{j-1}).$
	 Since
	$j-1\geq i$, the transitive action of $\Urt(n+1-2i)$ on $S^{2n+1-4i}$ implies the existence of $ g_i^1\in \Urt(n+1-2i)$ such that
	$g_i^1(0_{j-i-1},0,|\eta_1|,\eta_2,|\eta_3|,0,0_{j-i-1})=(1,0,0_{n-2i}).$
	Therefore, if $g_i=I_{\mathbb C^i}\times  g_i^1\times
	I_{\mathbb C^i}\in G_i$ then
	$g_ig_j\eta=(0_i,1,0,0_{n-i})\in S^{2n+1}.$
	
By repeating the same procedure for another element $\tilde \eta\in
	{S}^{2n+1}$, there exists $\tilde g_i\in G_i$ and $\tilde
	g_j\in G_j$ such that  $\tilde g_i\tilde g_j\tilde
	\eta=(0_i,1,0,0_{n-i})\in S^{2n+1}.$ Accordingly, $$\eta=g_j^{-1}g_i^{-1}\tilde
	g_i\tilde g_j \tilde \eta=g_j^{-1}\overline  g_i\tilde
	g_j \tilde \eta,$$ where $\overline g_i=g_i^{-1}\tilde
	g_i\in G_i,$ which concludes the proof. 
	\hfill $\square$\\
	
	For every fixed $i\in \{1,...,[\frac{n+1}{2}]\}$, we
	introduce the matrix $A_i$  as
		$$A_{i}=\left\{
	\begin{array}{lcl}
	\left[
	\begin{array}{cc}
	0 & I_{\mathbb C^{\frac{n+1}{2}}}  \\
	I_{\mathbb C^{\frac{n+1}{2}}} & 0  \\
	\end{array}
	\right],& \mbox{if}&
	n+1=2i,\\
	\\
	\left[
	\begin{array}{ccc}
	0 & 0 & I_{\mathbb C^{i}}  \\
	0& I_{\mathbb C^{n+1-2i}} & 0 \\
	I_{\mathbb C^{i}}& 0 & 0
	\end{array}
	\right],& \mbox{if}& n+1\neq 2i.
	\end{array}\right.$$
%	
%	$$\zeta_i=\left[
%	\begin{array}{cc}
%	0 & I_{\mathbb C^{\frac{n+1}{2}}}  \\
%	I_{\mathbb C^{\frac{n+2}{2}}} & 0  \\
%	\end{array}
%	\right]\  {\rm if}\  n+1=2i,  \ {\rm and}\  \zeta_i=\left[
%	\begin{array}{ccc}
%	0 & 0 & I_{\mathbb C^{i}}  \\
%	0& I_{\mathbb C^{n+1-2i}} & 0 \\
%	I_{\mathbb C^{i}}& 0 & 0
%	\end{array}
%	\right]\  {\rm if}\  n+1\neq 2i.
%	$$
%	\newpage
	The following construction is inspired  by  Bartsch and Willem
	\cite{BW-JFA}.  Since one has $A_i\in \Urt(n+1)\setminus G_i,$
	$A_i^2=I_{\mathbb C^{n+1}}$ and $A_i G_i=
	G_iA_i$, the group generated by
	$G_i$ and $A_i$ is $\hat G_i
	=[G_i; A_i]  =G_i\cup  A_iG_i$, i.e., 
{\small 	\begin{equation}\label{cseles-csoport}
	\hat G_i
	=\left\{
	\begin{array}{lcl}
	\left[
	\begin{array}{cc}
	{\bf U}(\frac{n+1}{2}) & 0  \\
	0 & {\bf U}(\frac{n+1}{2})  \\
	\end{array}
	\right]\cup \left[
	\begin{array}{cc}
	0 & 	{\bf U}(\frac{n+1}{2})  \\
		{\bf U}(\frac{n+1}{2}) & 0  \\
	\end{array}
	\right],& \mbox{if}&
	n+1=2i,\\
	\\
	\left[
	\begin{array}{ccc}
	{\bf U}(i) & 0 & 0  \\
	0& {\bf U}(n+1-2i) & 0 \\
	0& 0 & {\bf U}(i)
	\end{array}
	\right]\cup \left[
	\begin{array}{ccc}
	0 & 0 & {\bf U}(i)  \\
	0& {\bf U}(n+1-2i) & 0 \\
	{\bf U}(i)& 0 & 0
	\end{array}
	\right],& \mbox{if}& n+1\neq 2i.
	\end{array}\right.
	\end{equation}}
	In fact, in $\hat G_i$  there are only two types of elements: either of  the form $g\in G_i,$ or $A_ig \in 
	\hat G_i\setminus  G_i$ (with $g\in G_i$), respectively. 
	
	The action 
	$\hat G_i\circledast H^\gamma(S^{2n+1})\mapsto H^\gamma(S^{2n+1})$
	of the group $\hat G_i$ on
	$H^\gamma(S^{2n+1})$ is defined by
	\begin{equation}\label{zeta-action}
	(\hat g \circledast U)(\eta)=\left\{
	\begin{array}{ccl}
	U(g^{-1}\eta), & {\rm if} & \hat g=g\in G_i,\\
	-U(g^{-1}A_i^{-1}\eta), & {\rm if} & \hat g=A_ig \in \hat G_i\setminus  G_i,\\
	\end{array}
	\right.
	\end{equation}
	for every $\hat g\in \hat
	G_i,$ $U\in H^\gamma(S^{2n+1})$ and $\eta\in
	S^{2n+1}.$ 
	We notice that this action is well-defined, continuous and linear. Similarly as in (\ref{elso-action}), we introduce 
	the space of $G_i$-invariant functions of $H^\gamma(S^{2n+1})$ as
	$$H^\gamma_{G_i}(S^{2n+1})=\left\{U\in H^\gamma(S^{2n+1}): g\circ U=U\ \ \ {\rm for\ every}\ g\in G_i \right\},$$
	where the action $'\circ'$ corresponds to the first relation in (\ref{zeta-action}). Furthermore, let  
	$$H^\gamma_{\hat G_i}(S^{2n+1})=\left\{U\in H^\gamma(S^{2n+1}): \hat g \circledast U=U\ \ \ {\rm for\ every}\ \hat g\in \hat G_i \right\}$$
	be the space of $\hat  G_i$-invariant functions of $H^\gamma(S^{2n+1})$. 
	
	The following result summarizes the constructions in this subsection. 
	
	\begin{proposition}\label{proposition-Rubik} Let $\gamma\in \bigcup_{k=1}^n[k,\frac{kQ}{Q-1})$, and fix  $i, j\in
		\{1,...,[\frac{n+1}{2}]\}$ such that $i\neq j$. The following statements hold$:$
		\begin{enumerate}
			\item[{\rm (i)}] The embedding  $H^\gamma_{\hat G_i}(S^{2n+1})\hookrightarrow L^\frac{2Q}{Q-2\gamma}(S^{2n+1})$ is compact$;$
			
			\item[{\rm (ii)}]
			$H^\gamma_{G_i}(S^{2n+1})\cap H^\gamma_{G_j}(S^{2n+1})=\{constant\ functions\ on\ S^{2n+1}\};$ 
			
			\item[{\rm (iii)}] $H^\gamma_{\hat G_i}(S^{2n+1})\cap H^\gamma_{\hat G_j}(S^{2n+1})=\{0\}$.
			
			%\item[{\rm (ii)}] H^1_{T_{n,i}}({\mathbb R}^n)\cap
			%H^1_{T_{n,j}}({\mathbb R}^n)=H_{{\bf O}(n)}^1(\mathbb R^n)$
			%$H^1_{T^{\zeta_i}_{n,i}}({\mathbb R}^n)\cap
			%H^1_{T_{n,j}}({\mathbb R}^n)=\{0\}$ for every $i,j\in \{1,...,t_n\}$
			%with $i\neq j$;
		\end{enumerate}
	\end{proposition}
	{\it Proof.} (i) It is clear that $H^\gamma_{\hat G_i}(S^{2n+1})\subset H^\gamma_{G_i}(S^{2n+1})$.  Moreover, by Proposition \ref{proposition-compact} we have that
	the embedding $H_{G_i}^\gamma(S^{2n+1})\hookrightarrow L^\frac{2Q}{Q-2\gamma}(S^{2n+1})$ is compact.

	(ii) Let us fix $U\in H^\gamma_{G_i}(S^{2n+1})\cap H^\gamma_{G_j}(S^{2n+1}).$ Since $U$ is both $G_i$- and $G_j$- invariant, it is also $[G_i;G_j]$-invariant, i.e., $U(g\eta)=U(\eta)$ for every $g\in [G_i;G_j]$ and $\eta\in S^{2n+1}.$ According to Lemma \ref{Rubik}, the group $[G_i;G_j]$ acts transitively on the CR sphere
	${S}^{2n+1}$, i.e., the orbit of every element $\eta\in S^{2n+1}$ by the group $[G_i;G_j]$ is the whole sphere $S^{2n+1}$. Thus, $U$ is a constant function.

	(iii) Let $U\in H^\gamma_{\hat G_i}(S^{2n+1})\cap H^\gamma_{\hat G_j}(S^{2n+1}).$ On one hand, by (ii), we first have that $U$ is constant. On the other hand, the second relation from
	(\ref{zeta-action}) implies  that $U(\eta)=-U(A_i \eta)$ for every $\eta\in
	S^{2n+1}.$ Therefore, we necessarily have
	that $U=0$. \hfill $\square$\\

\medskip
	
	\subsection{Proof of Theorem \ref{fotetel}.}\label{section-3-4} 
	We associate to problem $(\textbf{FYS})_\gamma$ the energy functional $E:H^\gamma(S^{2n+1})\to \mathbb R$ defined by
	$$E(U)=\frac{1}{2}\int_{S^{2n+1}}U\mathcal A_\gamma U{\rm d}\eta-\frac{Q-2\gamma}{2Q}\int_{S^{2n+1}}|U|^\frac{2Q}{Q-2\gamma}{\rm d}\eta,\ \ \ \ U\in H^\gamma(S^{2n+1}).$$
	Due to (\ref{sharp-FL}), the functional $E$ is well-defined, belonging to $C^1(H^\gamma(S^{2n+1}), \mathbb R)$. Moreover,  $U\in H^\gamma(S^{2n+1})$ is a critical point of $E$ if and only if $U$ is a weak solution of $(\textbf{FYS})_\gamma$.

Let us fix $i\in
\{1,...,[\frac{n+1}{2}]\}$. In order to guarantee critical points for $E$, we first consider the functional $E_i:H^\gamma_{\hat G_i}(S^{2n+1})\to \mathbb R$, the restriction of $E$ to the space  $H^\gamma_{\hat G_i}(S^{2n+1})$. It is clear that $E_i$ is an even functional and it has the mountain pass geometry. Since the embedding  $H^\gamma_{\hat G_i}(S^{2n+1})\hookrightarrow L^\frac{2Q}{Q-2\gamma}(S^{2n+1})$ is compact, see Proposition \ref{proposition-Rubik} (i), we may apply the fountain theorem, see e.g. Bartsch and Willem \cite[Theorem 3.1]{BW-JFA}, guaranteeing a sequence $\{U_i^k\}_{k\in \mathbb N}\subset H^\gamma_{\hat G_i}(S^{2n+1})$  of critical points for $E_i$ with the additional property that $\|U_i^k\|_{H^\gamma}\to \infty$ as $k\to \infty$. 

By using the principle of symmetric criticality of Palais \cite{Palais}, we are going to prove that $\{U_i^k\}_{k\in \mathbb N}\subset H^\gamma_{\hat G_i}(S^{2n+1})$ are in fact critical points for the original energy functional $E$, thus  weak solutions of $(\textbf{FYS})_\gamma$. To do this, it suffices to verify that $E$ is a $\hat G_i$-invariant functional, i.e., $$E(\hat g \circledast U)=E(U) \ \ {\rm  for\ every}\ \hat g\in \hat G_i,\ U\in H^\gamma(S^{2n+1}).$$ 

 On one hand, according to relation (\ref{osszefugges}), for the quadratic  term in $E$  it is enough to prove that $\hat G_i$ acts isometrically on $H^\gamma(S^{2n+1})$, i.e., 
\begin{equation}\label{izometrikusan-hat}
(\hat g \circledast U,\hat g \circledast U)_\gamma=( U, U)_\gamma\ \ {\rm
for\ every}\ \hat g\in \hat G_i,\ U\in H^\gamma(S^{2n+1}).
\end{equation}
  To see this,   let us fix $\hat g\in \hat G_i$  and $U\in H^\gamma(S^{2n+1})$ arbitrarily. We recall that by definition
\begin{eqnarray*}
(\hat g \circledast U,\hat g \circledast U)_\gamma&=&\sum_{j,k\geq  0}\lambda_j(\gamma)\lambda_k(\gamma)\sum_{l=1}^{m_{j,k}}|c_{j,k}^l(\hat g \circledast U)|^2.
\end{eqnarray*}
By using (\ref{zonal-0}), one has 
\begin{eqnarray}\label{uros}
\nonumber\sum_{l=1}^{m_{j,k}}|c_{j,k}^l(\hat g \circledast U)|^2&=&\displaystyle\int_{S^{2n+1}}\int_{S^{2n+1}}(\hat g \circledast U)(\zeta)(\hat g \circledast U)(\eta)\sum_{l=1}^{m_{j,k}}Y_{j,k}^l(\zeta)\overline{Y_{j,k}^l(\eta)}{\rm d}\zeta{\rm d}\eta\\&=&\displaystyle\int_{S^{2n+1}}\int_{S^{2n+1}}(\hat g \circledast U)(\zeta)(\hat g \circledast U)(\eta)\Phi_{j,k}(\zeta,\eta){\rm d}\zeta{\rm d}\eta.
\end{eqnarray}
Note that for every $g\in G_i\subset {\bf U}(n+1)$ and $\zeta,\eta\in S^{2n+1}$ we have $$\langle g\zeta,\overline{g\eta}\rangle= \langle A_ig\zeta,\overline{A_ig\eta}\rangle=\langle \zeta,\overline{\eta}\rangle;$$ therefore, by 
the representation (\ref{zonal}) of the zonal harmonics we also have that
$$\Phi_{j,k}(g\zeta,{g\eta})=\Phi_{j,k}(A_ig\zeta,{A_ig\eta})=\Phi_{j,k}(\zeta,\eta).$$
 Thus, relation (\ref{zeta-action}) and suitable changes of variables in (\ref{uros}) imply that $$\sum_{l=1}^{m_{j,k}}|c_{j,k}^l(\hat g \circledast U)|^2=\displaystyle\int_{S^{2n+1}}\int_{S^{2n+1}}U(\zeta)U(\eta)\Phi_{j,k}(\zeta,\eta){\rm d}\zeta{\rm d}\eta=\sum_{l=1}^{m_{j,k}}|c_{j,k}^l(U)|^2,$$ 
 which proves (\ref{izometrikusan-hat}). 
 
 On the other hand, the $\hat G_i$-invariance of the nonlinear term $U\mapsto \displaystyle\int_{S^{2n+1}}|U|^\frac{2Q}{Q-2\gamma}$ trivially follows by a change of variable, by using the isometric character of the group ${\bf U}(n+1)$ on $S^{2n+1}$.
 
 Accordingly, for every $i\in
 \{1,...,[\frac{n+1}{2}]\}$, the functions $\{U_i^k\}_{k\in \mathbb N}\subset H^\gamma_{\hat G_i}(S^{2n+1})$ are non-trivial weak solutions of $(\textbf{FYS})_\gamma$. Due to Proposition  \ref{weak-equivalent}, 
 $u_i^k=(2{\rm Jac}_{\mathcal C})^\frac{Q-2\gamma}{2Q}U_i^k\circ \mathcal C\in  D^\gamma(\mathbb H^n)$
 are non-trivial weak solutions of the original fractional Yamabe problem {\rm${(\textbf{FYH})_\gamma}$}; by construction, $u_i^k$ are sign-changing functions.   

Due to Proposition \ref{proposition-Rubik} (iii), we state that the sequences $\{U_i^k\}_{k\in \mathbb N}\subset H^\gamma_{\hat G_i}(S^{2n+1})$ and $\{U_j^k\}_{k\in \mathbb N}\subset H^\gamma_{\hat G_j}(S^{2n+1})$ with $i,j\in
\{1,...,[\frac{n+1}{2}]\}$, $i\neq j,$ cannot be compared from symmetrical point of view. Therefore, the sequences $\{u_i^k\}\subset D^\gamma(\mathbb H^n)$ and $\{u_j^k\}\subset D^\gamma(\mathbb H^n)$ have mutually different nodal properties for every $i,j\in
\{1,...,[\frac{n+1}{2}]\}$, $i\neq j,$ which concludes the proof.

	\begin{remark}\rm \label{rem-sign-chang-1}
Consider  a nonzero solution  $u_i^k=(2{\rm Jac}_{\mathcal C})^\frac{Q-2\gamma}{2Q}U_i^k\circ \mathcal C\in  D^\gamma(\mathbb H^n)$ of $(\textbf{FYH})_\gamma$, with  $\{U_i^k\}_{k\in \mathbb N}\subset H^\gamma_{\hat G_i}(S^{2n+1})\setminus \{0\}$. For simplicity, we consider the case $n+1=2i$.  Let us introduce the nodal domain of $U_i^k$ (or $u_i^k$) as the connected components of $C_i^k=S^{2n+1}\setminus N_i^k$, where  $N_i^k=\overline{\{\eta\in S^{2n+1}:U_i^k(\eta)=0\}}.$ Since $U_i^k\in  H^\gamma_{\hat G_i}(S^{2n+1}),$ by relation (\ref{zeta-action}) it follows that $U_i^k$ has the form $U_i^k(\eta)=U_i^k(|\eta_1|,|\eta_2|)$ with the property that $U_i^k(|\eta_1|,|\eta_2|)=-U_i^k(|\eta_2|,|\eta_1|)$, $\eta=(\eta_1,\eta_2)\in S^{2n+1}$, $\eta_1,\eta_2\in \mathbb C^i.$ Accordingly, since $U_i^k(|\eta_1|,|\eta_2|)=U_i^k(|\pm\eta_1|,|\pm\eta_2|)$, $U_i^k$ is sign-changing with at least four non-degenerate nodal domains in $C_i^k$; in two of them the function $U_i^k$ is negative, while in the other two it is positive, respectively. When $n+1\neq 2i,$ a similar discussion can be performed.  
\end{remark}

We conclude the paper by the following table providing explicit forms of subgroups of the unitary group ${\bf U}(n+1)$ and admissible intervals for the parameter $\gamma$, depending on the dimension $n$, where our main theorem applies; we only consider the cases when $n\in \{1,...,8\}$:

{\footnotesize
	\begin{center}
		\begin{tabular}[b]{|c|c|c|c|c|}
			\hline
			\cline{1-5} $n$ &$Q=2n+2$&  $G_i$, $i\in \{1,...,[\frac{n+1}{2}]\}$ & Admissible domains & Number of symmetrically distinct\\
			&&&for $\gamma\in (0,Q/2)$& sequences of solutions of {\rm${(\textbf{FYH})_\gamma}$}
			\\ \hline\cline{1-5}
			
			$1$&$4$&$G_{1}={\bf U}(1)\times {\bf U}(1)$& $[1,\frac{4}{3})$ & 1\\
			
			 \cline{1-5}
			
			$2$&$6$&$G_{1}=\Urt(1)\times \Urt(1)\times
			\Urt(1)$&$[1,\frac{6}{5})\cup [2,\frac{12}{5})$&1\\
			
			\cline{1-5}
			
			$3$&$8$&$G_{1}=\Urt(1)\times \Urt(2)\times
			\Urt(1)$&$[1,\frac{8}{7})\cup [2,\frac{16}{7})\cup $&2\\
			
			$$&$$&$G_{2}=\Urt(2)\times \Urt(2)$&$\cup [3,\frac{24}{7})$ &\\

			\cline{1-5}
			
			$4$&$10$&$G_{1}=\Urt(1)\times \Urt(3)\times
			\Urt(1)$&$\bigcup_{k=1}^4[k,\frac{10k}{9})$&2\\
			
			$$&$$&$G_{2}=\Urt(2)\times \Urt(1)\times
			\Urt(2)$&&\\
			
			\cline{1-5}
			
			$$&$$&$G_{1}=\Urt(1)\times \Urt(4)\times
			\Urt(1)$&&\\
			
			$5$&$12$&$G_{2}=\Urt(2)\times \Urt(2)\times
			\Urt(2)$&$\bigcup_{k=1}^5[k,\frac{12k}{11})$&3\\
			
			$$&$$&$G_{3}=\Urt(3)\times \Urt(3)$&&\\

			\cline{1-5}
			
			$$&$$&$G_{1}=\Urt(1)\times \Urt(5)\times
			\Urt(1)$&&\\
			
			$6$&$14$&$G_{2}=\Urt(2)\times \Urt(3)\times
			\Urt(2)$&$\bigcup_{k=1}^6[k,\frac{14k}{13})$&3\\
			
			$$&$$&$G_{3}=\Urt(3)\times \Urt(1)\times
			\Urt(3)$&&\\
			
			\cline{1-5}
			
			$$&$$&$G_{1}=\Urt(1)\times \Urt(6)\times
			\Urt(1)$&&\\
			
			$7$&$16$&$G_{2}=\Urt(2)\times \Urt(4)\times
			\Urt(2)$&$\bigcup_{k=1}^7[k,\frac{16k}{15})$&4\\
			
			$$&$$&$G_{3}=\Urt(3)\times \Urt(2)\times
			\Urt(3)$&&\\
			
			$$&$$&$G_{4}=\Urt(4)\times \Urt(4)$&&\\

			\cline{1-5}
			
			$$&$$&$G_{1}=\Urt(1)\times \Urt(7)\times
			\Urt(1)$&&\\
			
			$8$&$18$&$G_{2}=\Urt(2)\times \Urt(5)\times
			\Urt(2)$&$\bigcup_{k=1}^8[k,\frac{18k}{17})$&4\\
			
			$$&$$&$G_{3}=\Urt(3)\times \Urt(3)\times
			\Urt(3)$&&\\
			
			$$&$$&$G_{4}=\Urt(4)\times \Urt(1)\times
			\Urt(4)$&&\\
			\cline{1-5}
	\end{tabular}\end{center}}

	\vspace*{0.6cm}
	\noindent {\bf Acknowledgment.} The author would
	like to thank Professor Carlo Morpurgo  for his useful remarks concerning the preliminary version of the manuscript and the  anonymous referee for her/his valuable observations which improved the quality of the manuscript.

\end{document}